\documentclass[12pt]{article}
\usepackage{amsfonts,amsmath,amssymb,amscd,theorem}
\newtheorem{definition}{Definition}[section]

\newtheorem{proposition}[definition]{Proposition}
\newtheorem{theorem}[definition]{Theorem}

\theorembodyfont{\rmfamily}
\newtheorem{rmk}{Remark}[section]

\def\cA{{\cal A}}          \def\cB{{\cal B}}          
          \def\cE{{\cal E}}          \def\cF{{\cal F}}
          \def\cH{{\cal H}}

                    \def\cR{{\cal R}}
                    \def\cU{{\cal U}}

\def\fa{{\mathfrak a}}

\newcommand{\CC}{{\mathbb C}}
\newcommand{\II}{{\mathbb I}}
\newcommand{\NN}{{\mathbb N}}
\newcommand{\ZZ}{{\mathbb Z}}
\newcommand{\eps}{{\varepsilon}}

\newcommand{\un}{{\mathbb I}}
\newcommand{\id}{1}

\newcommand{\sfrac}[2]{{\textstyle{\frac{#1}{#2}}}}
\newcommand{\half}{{\sfrac{1}{2}}}
\newcommand{\three}{{\sfrac{3}{2}}}

\newcommand{\finproof}{{\hfill \rule{5pt}{5pt}}}

\def\QTHA{Quasitriangular Hopf Algebra (QTHA)\def\QTHA{QTHA}}
\def\QTQHA{Quasitriangular Quasi-Hopf Algebra (QTQHA)\def\QTQHA{QTQHA}}

\begin{document}
\pagestyle{empty}

\markright{\today\dotfill DRAFT\dotfill }

\begin{center}

{\Large \textsf{
Super Yangian $Y(osp(1|2))$ \\[4mm] 
and the \\[5mm]
Universal $R$-matrix of its Quantum Double
}}

\vspace{10mm}

{\large D. Arnaudon, N.~Cramp\'e, L. Frappat$^{*}$, {E}.
Ragoucy}

\vspace{10mm}

\emph{Laboratoire d'Annecy-le-Vieux de Physique Th{\'e}orique}

\emph{LAPTH, CNRS, UMR 5108, Universit{\'e} de Savoie}

\emph{B.P. 110, F-74941 Annecy-le-Vieux Cedex, France}

\vspace{7mm}

\emph{$^*$ Member of Institut Universitaire de France}

\end{center}

\vfill
\vfill

\begin{abstract}
We present the Drinfel'd realisation of the super Yangian
$Y(osp(1|2))$, including the explicit expression for the coproduct. We
show in particular that it is necessary to introduce 
supplementary Serre relations. The construction of its quantum double
is carried out. This allows us to give the
universal $R$-matrix of $DY(osp(1|2))$. 
\end{abstract}

\vfill
\vfill

\rightline{LAPTH-932/02}
\rightline{math.QA/0209167}
\rightline{September 2002}

\baselineskip=16pt


\newpage

\pagestyle{plain}
\setcounter{page}{1}

\section{Introduction}

The Yangian $Y(\fa)$ based on a simple Lie algebra $\fa$ is defined as the
homogeneous quantisation of the algebra $\fa[u] = \fa \otimes
\CC[u]$ endowed with its standard bialgebra structure, where $\CC[u]$ is
the ring of polynomials in the indeterminate $u$. It was introduced by
Drinfel'd \cite{Dri85,Dri86,Dri88}. The most elegant and concise
presentation of the Yangian uses the FRT formalism \cite{FRT}, based
on a certain 
evaluated R-matrix. For unitary algebras, the R-matrix is given by $R(u) =
\II\otimes\II + P/u$ where $P$ is the permutation map. This can be extended to
superunitary series by considering the superpermutation instead
\cite{KS,Naz91}. This 
formalism was extended to the orthogonal, symplectic and orthosymplectic
cases by taking $R(u) = \II\otimes\II + P/u - K/(u+\kappa)$ where $K$
is a partial 
(super)transposition of $P$ and $4\kappa = (\alpha_0+2\rho,\alpha_0)$,
see for example 
\cite{soya}. Although one can exhibit a suitable R-matrix for any simple
Lie algebra or basic simple Lie superalgebra $\fa$ for defining the Yangian
$Y(\fa)$, the obtained structure is a Hopf algebra, but not a
quasi-triangular one. It is well known that the quantum double construction
allows one to construct a quasi-triangular Hopf algebra from a Hopf
algebra, and that this procedure leads to the universal R-matrix of the
algebra under consideration. However, in order to find explicit useful
formula for this universal R-matrix, it is necessary to consider another
realisation of the Yangian, given in terms of generators and relations
similar to the description of a loop algebra as a space of maps.
Unfortunately in this realisation no explicit formula for the
comultiplication is known in general, except in the $sl(2)$ case
\cite{Molev2001}.

The aim of this paper is to extend this construction to the case of the
superalgebra $osp(1|2)$. Let us recall that the super Yangian $Y(osp(1|2))$
is defined by the relations
\begin{eqnarray}
&& R_{12}\left(u-v\right) \, L_1(u) \, L_2(v) = L_2(v) \, L_1(u) \,
R_{12}\left(u-v\right) \\
&& \hspace*{-48mm} \mbox{and} \nonumber \\
&& C(u) = L^t\left(u-\kappa\right) \, L(u)=\id_{3}
\end{eqnarray}
where the generators of $Y(osp(1|2))$ are encapsulated into the $3
\!\times\! 3$ matrix $L(u)$ and $R_{12}(u)$ is the R-matrix introduced in
\cite{soya}. Writing a Gauss decomposition of $L(u)$, we then introduce
some specific combinations of the generators $L^{ij}(u)$, called $e(u)$,
$f(u)$ and $h(u)$, which define a Drinfel'd realisation of the
super Yangian $Y(osp(1|2))$. More precisely, we show that the associative
algebra $\cA$ generated by $e(u)$, $f(u)$ and $h(u)$ subjected to certain
relations, and the super Yangian $Y(osp(1|2))$ are isomorphic as
bialgebras. At this point, two remarks are in order. First, one is able to
find explicit formula for the comultiplication in terms of the Drinfel'd
generators $e(u)$, $f(u)$ and $h(u)$, thereby generalising Molev's formula
in the case of $osp(1|2)$. Second, it is necessary to introduce
supplementary Serre-type relations among the $e(u)$, $f(u)$, $h(u)$
generators, as in the case of $U_{q}(\widehat{osp}(1|2))$ \cite{Ding}.
These supplementary Serre-type relations are cubic in the Drinfel'd
generators. Indeed, it appears that the quadratic exchange relations among
the $e(u)$, $f(u)$, $h(u)$ generators, derived from the RLL relations, lead
to a superalgebra which is bigger than $Y(osp(1|2))$. Hence it is necessary
to quotient this bigger structure by supplementary relations.

The next step is the construction of the quantum double. The super Yangian
$Y(osp(1|2)) = Y^{+}$ being given in terms of the Drinfel'd generators
$e(u)$, $f(u)$, $h(u)$ with positive modes, we introduce another set of
generators $e(u)$, $f(u)$, $h(u)$ with negative modes, generating a Hopf
algebra $Y^{-}$. We construct a Hopf pairing between $Y^{+}$ and $Y^{-}$,
such that $Y^{-}$ is isomorphic to the dual of $Y^{+}$ with opposite
comultiplication. We prove that this Hopf pairing is not degenerate. This
allows us to define the double super Yangian $DY(osp(1|2))$ in a proper
way. In particular we are able to give a presentation of the double
super Yangian $DY(osp(1|2))$ in terms of the Drinfel'd generators $e(u)$,
$f(u)$, $h(u)$, now $\ZZ$ moded, subjected to suitable quadratic exchange
relations, and suitable supplementary Serre-type relations of cubic form.

As we emphasised above, the quantum double procedure allows one to
construct explicitly the universal R-matrix. $\{x_i,i\in\NN\}$ being a
basis of $Y^{+}$ and $\{x^i,i\in\NN\}\in Y^{-}$ its dual basis (i.e.
$<x^i,x_j>=\delta^i_j$), the universal R-matrix of $DY$ is given by
$\cR=\sum x_i \otimes x^i$. Let $\widetilde\cE^+$, $\widetilde\cF^+$ and
$\widetilde\cH^+$ denote the unital subalgebras of $Y^+$ generated by the
positive modes of $e(u)$, $f(u)$ and $h(u)$, and let $\widetilde\cF^-$,
$\widetilde\cE^-$ and $\widetilde\cH^-$ be the dual subalgebras. Following
the kind of arguments used in \cite{Rosso}, we give
Poincar\'e--Birkhoff--Witt bases for $\widetilde\cE^\pm$,
$\widetilde\cF^\pm$ in terms of the modes of the generators $e(u)$ and
$f(u)$. We show that this leads to the usual nice factorised expression of
the universal R-matrix of $DY(osp(1|2))$, namely $\cR=\cR_E\cR_H\cR_F$,
where $\cR_E\in \widetilde\cE^+ \otimes \widetilde\cF^-$, $\cR_H \in
\widetilde\cH^+ \otimes \widetilde\cH^-$ and $\cR_F \in
\widetilde\cF^+\otimes\widetilde\cE^-$. Finally, considering the action of
the universal R-matrix $\cR$ on evaluation representations of
$DY(osp(1|2))$, we obtain the evaluated R-matrix of $DY(osp(1|2))$, which
coincides with the one introduced in \cite{soya}, up to a normalisation
factor written as a ratio of $\Gamma_{1}$ functions of period $2\kappa$.

\section{The RTT presentation of super Yangian $Y(osp(1|2))$}
\setcounter{equation}{0}

We denote by $V$ the 3-dimensional $\ZZ_2$-graded vector space
representation of $osp(1|2)$. The first and the third  
basis vectors have the grade 0 (mod 2) whereas the second has the
grade 1 (mod 2). The same gradation was used  
by Ding \cite{Ding} in order to define $U_q(\widehat{osp}(1|2))$. The
multiplication for the tensor product is defined for $a$, $b$,  
$\alpha$, $\beta \in Y(osp(1|2))$ by
\begin{equation}(a \otimes \alpha)(b \otimes \beta)=(-1)^{[b][\alpha]} 
(a b\otimes \alpha \beta)
\end{equation} 
where $[a]\in \ZZ_2$ denotes the grade of $a$. Let $E_{ij}$ be the  
elementary matrix with entry 1 in row $i$ and column $j$ and 0
elsewhere. The ``usual'' super transposition $.^T$ is defined by
$A^T=\sum_{i,j=1}^3 (-1)^{[j]([i]+[j])}A^{ji}E_{ij}$ for any matrix 
$A=\sum_{i,j=1}^3 A^{ij}E_{ij}$.
The super transposition $.^t$ we will use is a conjugation of the previous one:
\begin{equation}
A^t=\sum_{i,j,k,l=1}^3(-1)^{[i][l]+[i]}J^{ij}A^{kj}J^{lk}E_{il}=J\,A^T\,J^{-1}\
\mbox{where}\; J=E_{31}+E_{22}-E_{13}
\end{equation}
The super permutation $P_{12}$ (i.e. $X_{21}=P_{12}X_{12}P_{12}$)
is defined
by $P_{12}=(-1)^{[j]}E_{ij}\otimes E_{ji}$.\\
The super R-matrix $R_{12}(u)\in End(V \otimes V)$ is defined by: 
 
\begin{eqnarray}
  R_{12}(u)&=&\rho(u) \frac{u}{u+1}\left({\un}\otimes {\un}
+\frac{P}{u}-\frac{P^{t_1}}{u+\kappa}\right)\\[7pt]
  &=&\rho(u) \frac{u}{u+1}
  \left(\begin{array}{c c c c c c c c c}
      \frac{u+1}{u} & 0 & 0& 0& 0& 0& 0& 0&0\\
      0& 1 & 0 &-\frac{1}{u}& 0 & 0 & 0 & 0 &0\\
      0& 0 & \frac{u+\kappa-1}{u+\kappa} & 0 &\frac{1}{u+\kappa}
      & 0 & \frac{\left(2u+\kappa\right)}{u(u+\kappa)}& 0 &0\\ 
      0&\frac{1}{u}& 0 & 1 & 0 & 0 & 0 & 0& 0\\
      0& 0 &\frac{1}{u+\kappa} &  
      0&\frac{u^2+\kappa u-\kappa}{u(u+\kappa)} 
      & 0 &-\frac{1}{u+\kappa}
      & 0 &0\\ 
      0& 0 & 0 & 0 & 0 & 1 & 0 &\frac{1}{u}&0\\
      0& 0 &\frac{\left(2u+\kappa\right)}{u(u+\kappa)}
      & 0 &-\frac{1}{u+\kappa}& 0 &\frac{u+\kappa-1}{u+\kappa} & 0 &0\\
      0& 0 & 0 & 0 & 0 &-\frac{1}{u}& 0 & 1 &0\\
      0&0&0&0&0&0&0&0&\frac{u+1}{u}
    \end{array}
  \right) \label{matrixR} 
\end{eqnarray}
where  
\begin{equation}
\rho(u) = 
\frac{\Gamma_{1}\left(u|2\kappa\right)\,
    \Gamma_{1}\left(u+\kappa-1|2\kappa\right)\,
    \Gamma_{1}\left(u+\kappa+1|2\kappa\right)\,
\Gamma_{1}\left(u+2\kappa|2\kappa\right)} 
  {\Gamma_{1}\left(u+1|2\kappa\right)\,
    \Gamma_{1}\left(u+\kappa|2\kappa\right)^{2}\,
    \Gamma_{1}\left(u+2\kappa-1|2\kappa\right)} \;, 
\;\; \kappa=\frac{3}{2}\;,
  \label{matrice_R}
\end{equation}
$.^{t_1}$ is the super transposition in the first space and 
$\un$ is the $3\times3$ identity matrix.
The function $\Gamma_1$ is defined by
$
  \Gamma_1(x|\omega) = 
  \frac{\displaystyle\omega^{x/\omega}}
  {\displaystyle \sqrt{2\pi \omega}}
  \Gamma\left(\frac{\displaystyle x}{\displaystyle\omega}\right)
$.

It is known (see for instance \cite{soya}) that:
\begin{proposition}
  The matrix $R(u)$ satisfies
  \begin{eqnarray}
    &&R_{12}(u) \, R_{13}(u+v) \, R_{23}(v) = R_{23}(v) \, R_{13}(u+v)
    \, R_{12}(u),\qquad\mbox{(super Yang--Baxter)} \label{YBE}\\
    &&R_{12}^{t_{1}}\left(-u-\kappa\right)=
    R_{12}(u),\qquad \mbox{(crossing symmetry)} \label{eq:propR} \\ 
    &&R_{12}(u) \, R_{12}(-u) = \rho(u)\rho(-u) , 
\qquad\mbox{(unitarity)}\label{eq:Runit}
  \end{eqnarray}
\end{proposition}

We gave in our previous paper \cite{soya} the RTT presentation of
super Yangian $Y(osp(1|2))$. In the following,
$\id$ will denote the unit of the algebra  and
$\id_3=\id\,\un$. 
\begin{theorem}
  The super Yangian $Y(osp(1|2))$ is isomorphic to the associative
  superalgebra $\cU(R)$ generated by the elements $L^{ij}_{(n)}$ 
  ( $1\le i,j\le 3$, $n\in\ZZ_{>0}$ ), $\id$ and the defining
  relations, given in terms of formal series 
  $L(u) 
  =\id_3 + \sum_{i,j=1}^{3} \sum_{n \in \ZZ_{> 0}}
  L^{ij}_{(n)} \, u^{-n} \, E_{ij}
  = \sum_{i,j=1}^{3} {L}^{ij}(u) \, E_{ij}$ :
  \begin{eqnarray}
    \bullet &&R_{12}\left(u-v\right) \, L_1(u) \, L_2(v)
    = L_2(v) \, L_1(u) \, R_{12}\left(u-v\right)\label{RLL}\\    
    \bullet &&C(u) = L^t\left(u-\kappa\right) \,
    L(u)=\id_{3}\label{C}
  \end{eqnarray}
  The Hopf algebra structure of $\cU(R)$ is given by
  \begin{eqnarray}
    \label{eq:coproduct}
    && \Delta \big(L(u) \big) = L(u) \,\dot\otimes\,   L(u) 
    \qquad \mbox{i.e.}\qquad
    \Delta \big( {L}^{ij}(u) \big) = \sum_{k=1}^{3} {L}^{ik}(u) \otimes
    {L}^{kj}(u) 
    \\
    && S(L(u)) =  L(u)^{-1}
    \quad ; \quad
    \varepsilon(L(u)) = \un
    \label{eq:antipode}
  \end{eqnarray}
\end{theorem}

\section{The Drinfel'd realisation of $Y(osp(1|2))$}
\label{sect:A+}
\setcounter{equation}{0}

\begin{definition}\label{A+alg}
Let $\cA^+$ be the associative superalgebra generated by the odd elements
e$_k$, f$_k$ ($k\in\ZZ_{\ge0}$), the even elements  h$_k$
($k\in\ZZ_{\ge0}$), the unit $\id$ and 
the defining relations, given in terms of the generating functions 
$
e(u)=\sum_{k=0}^{\infty}e_ku^{-k-1},\;
f(u)=\sum_{k=0}^{\infty}f_ku^{-k-1},\;
h(u)=\id+ \sum_{k=0}^{\infty}h_ku^{-k-1} : 
$
\begin{eqnarray}
  \label{generatrice} 
  [h(u),h(v)] & =& 0,\label{gen_hh}\\
  \{e(u),f(v)\}&=&\frac{h(v)-h(u)}{u-v},\\
  \,[h(u),e(v)]
  &=& -\frac{(u-v+\kappa-1)[h(u),e(u)-e(v)]}{(u-v)(u-v+\kappa)}\nonumber\\
  &&+\frac{h(u)(e(v)+e(u)-e_0)-(e(u)-e_0)h(u)}{u-v+\kappa}\\
  \,[h(u),f(v)]  
  &=& \frac{(u-v-\kappa+1)[h(u),f(u)-f(v)]}{(u-v)(
  u-v-\kappa)} \nonumber\\
  &&-\frac{h(u)(f(v)+f(u)+f_0)-(f(u)+f_0)h(u)}{u-v-\kappa}\\
  \{e(u),e(v)\}  &=&\frac{[e(u),e(v)]}{2\,(u-v)}-
  \frac{\{e_0,e(u)-e(v)\}}{(u-v)}- \frac{(e(u)-e(v))^2}{2\,(u-v)^2},\\ 
  \{f(u),f(v)\}  &=&-\frac{[f(u),f(v)]}{2\,(u-v)}-
  \frac{\{f_0,f(u)-f(v)\}}{(u-v)}- \frac{(f(u)-f(v))^2}{2\,(u-v)^2}  \label{gen_ff}
\end{eqnarray}
and the supplementary Serre relations 
\begin{eqnarray}
  \label{eq:serree}
  &&  e(u)^3  =   e(u)\{e(u),e_0\}+[e_0^2,e(u)] \;,\\  
  \label{eq:serref}
  &&  f(u)^3  =  - f(u)\{f(u),f_0\}+[f_0^2,f(u)] \;.
\end{eqnarray}
\end{definition}   
The relations (\ref{gen_hh})--(\ref{gen_ff}) are equivalent to the
following commutation relations in terms of the modes 
$e_k$, $f_k$, $h_k$ ($k\ge0$)  :
\begin{eqnarray}
\bullet\, h_k\; \mbox{and}\; h_l:\ \  && [h_k,h_l]=0 ,
\label{h_h}\\[1.2ex]
\bullet\, e_k\; \mbox{and}\; f_l:\ \ &&\{e_k,f_l\}-h_{k+l}=0 ,
\label{e_f}\\[1.2ex]
\bullet\,  h_k \; \mbox{and}\; e_l :\ \ &&[h_0,e_l]-e_l=0 ,
\label{h0_e}\\[1.2ex]
&& 2\,[h_1,e_l]-2 e_{l+1}=\{h_0,e_l\},
\label{h1_e}\\[1.2ex]
&&2\,[h_{k+2},e_l]+2\,[h_{k},e_{l+2}]-4\,[h_{k+1},e_{l+1}]
=   [h_k,e_l]+\{h_{k+1},e_l\}-\{h_k,e_{l+1}\},
\ \ \ \ \ 
\label{h_e} 
\end{eqnarray}
\begin{eqnarray}
\bullet\, h_k\; \mbox{and}\; f_l :\ \ && [h_0,f_l]+f_l=0 ,
\label{h0_f}\\[1.2ex]
    &&2\,[h_1,f_l]+2 f_{l+1}=-\{h_0,f_l\},
\label{h1_f}\\[1.2ex] 
    &&2\,[h_{k+2},f_l]+2\,[h_{k},f_{l+2}]-4\,[h_{k+1},f_{l+1}]
    =[h_k,f_l]-\{h_{k+1},f_l\}+
\{h_k,f_{l+1}\},
\label{h_f} \ \ \ \ \ 
\end{eqnarray}
\begin{eqnarray}
\bullet\, e_k \; \mbox{and}\; e_l :\ \ &&
2\{e_{k+2},e_l\}+2\{e_{k},e_{l+2}\}-4\{e_{k+1},e_{l+1}\}
=\{e_k,e_l\}+[e_{k+1},e_l]-[e_k,e_{l+1}],
\label{e_e}\ \ \ \ \ 
\end{eqnarray}
\begin{eqnarray}
\bullet\, f_k\; \mbox{and}\; f_l :\ \ &&
2\{f_{k+2},f_l\}+2\{f_{k},f_{l+2}\}-4\{f_{k+1},f_{l+1}\}
=\{f_k,f_l\}-[f_{k+1},f_l]+[f_k,f_{l+1}]
    .\label{f_f}\ \ \ \ \ \ 
\end{eqnarray}
The Serre relations (\ref{eq:serree}) and (\ref{eq:serref}) in terms
of modes are conjectured to be (for $k\ge 0$)
\begin{eqnarray}
    [\{e_k,e_{k+1}\},e_k] &=& 2 e_k^3,
\label{e_E1}\\[1.2ex]
[\{e_k,e_{k+1}\},e_{k+1}] &=& - e_k e_{k+1} e_k - \frac{1}{2} e_k^2
    e_{k+1} - \frac{1}{2}e_{k+1} e_k^2,
\label{e_E2}\\[1.2ex] 
    [\{e_k,e_{k+1}\},e_{k+2}] &=& -2 e_{k+1}^2 e_k-2 e_{k+1} e_k
    e_{k+1}-2 e_k e_{k+1}^2,
\label{e_E3} 
\end{eqnarray}
\begin{eqnarray}
  [\{f_k,f_{k+1}\},f_k] &=& -2 f_k^3,
\label{f_F1}\\[1.2ex]
[\{f_k,f_{k+1}\},f_{k+1}] &=&  f_k f_{k+1} f_k + \frac{1}{2} f_k^2
    f_{k+1} + \frac{1}{2}f_{k+1} f_k^2,
\label{f_F2}\\[1.2ex] 
[\{f_k,f_{k+1}\},f_{k+2}] &=& 2 f_{k+1}^2 f_k+2 f_{k+1} f_k
    f_{k+1}+2 f_k f_{k+1}^2 \;.
\label{f_F3} 
\end{eqnarray}
This conjecture is supported by two results: on the one hand we have
proved them in the graded algebra (to be defined below); on the other
hand we checked 
explicitly the first nine relations. 
\begin{proposition}
\label{grad}
The algebra, $\cA^+$ (resp. $\cU(R)$), can be equipped 
with an ascending filtration with 
the degree of the generators defined by $deg(e_k)=k$, $deg(f_k)=k$, 
$deg(h_k)=k$ ($deg(L^{ij}_{(k)})=k-1$) and $deg(xy)=deg(x)+deg(y)$, 
for $x,y \in \cA^+$ (resp. $x,y \in \cU(R)$). \\
Let $gr(\cA^+)$ and $gr(\cU(R))$ denote the corresponding graded algebras
and $osp(1|2)[u]$ denote the Lie super algebra of polynomials 
in an indeterminate u with coefficients in $osp(1|2)$.
The algebras $gr(\cA^+)$, $gr(\cU(R))$ and $U(osp(1|2)[u])$ are 
isomorphic.
\end{proposition}
\textbf{Proof:} We first recall the notion of graded algebras.
We start with an algebra $\cA$, equipped with a grading $deg$, 
i.e. a morphism from
$(\cA,.)$ to $(\NN,+)$.
One introduces $\cA_k=\{x\in\cA, deg(x)\leq k\},\ k\geq0$ and
$gr(\cA_k)=\cA_k/\cA_{k-1}$ for $k\geq1$, $gr(\cA_0)=\cA_0$.
Then the graded algebra of $\cA$ is 
$gr(\cA)=\oplus_{k\geq0}gr(\cA_k)$.\\
The algebra $gr(\cA^+)$ is the algebra generated by 
$e_k, f_k, h_k\,(k\in\ZZ)$ 
and the relations (\ref{h_h})-(\ref{f_F3}) where the right hand 
side of the equalities is substituted by zero. These equalities are 
equivalent to :
\begin{eqnarray}
  [h_k,h_l]=0 &,& \{e_k,f_l\}=h_{k+l}\;,\\
  {}[h_k,e_l]=e_{k+l}&,&[h_k,f_l]=-f_{k+l}\;,\\
  \{e_n,e_m\}=\{e_0,e_{n+m}\}&,&\{f_n,f_m\}=\{f_0,f_{n+m}\}\;,\\
  {}[\{e_l,e_m\},e_n]=0&,&[\{f_l,f_m\},f_n]=0\;.
\end{eqnarray}
which are the relations of $U(osp(1|2)[u])$. Then, $gr(\cA^+)$ is
isomorphic to $U(osp(1|2)[u])$.\\
The isomorphism between 
$gr(\cU(R))$ and $U(osp(1|2)[u])$ is proved in \cite{soya}.
\finproof
\begin{theorem}
\label{iso}
The linear map
\begin{eqnarray}
  \phi:\hspace{20pt}\cA^+&\longrightarrow& \cU(R) \nonumber\\
  e(-u)&\longmapsto&{L}^{33}(u)^{-1}{L}^{23}(u) \label{isoe}\\
  f(-u)&\longmapsto&{L}^{32}(u){L}^{33}(u)^{-1} \label{isof}\\
  h(-u)&\longmapsto&{L}^{22}(u){L}^{33}(u)^{-1}+{L}^{32}(u)\label{isoh}
  {L}^{33}(u)^{-1}{L}^{23}(u){L}^{33}(u)^{-1}
\end{eqnarray}
is an isomorphism of algebra.\\
\end{theorem}
\textbf{Proof:}
The first step of the proof is to show that $\phi$ is a morphism of algebra.
\begin{eqnarray}
  &&\{\phi(e(-u)),\phi(f(-v))\}=\{{L}^{33}(u)^{-1}{L}^{23}(u),
  {L}^{32}(v){L}^{33}(v)^{-1}\} \nonumber\\
  &=&-{L}^{33}(u)^{-1}{L}^{32}(v)\left[{L}^{23}(u),{L}^{33}(v)^{-1}\right]
  +{L}^{33}(u)^{-1}\left\{{L}^{23}(u),{L}^{32}(v)\right\}{L}^{33}(v)^{-1}
   \nonumber\\
  &&-\left[{L}^{33}(u)^{-1},{L}^{32}(v)\right]{L}^{33}(v)^{-1}{L}^{23}(u)
   \\
  &=&\frac{1}{u-v}\left(-{L}^{33}(u)^{-1}\phi(h(-v))
    {L}^{33}(u)+{L}^{33}(u)^{-1}\phi(h(-u)){L}^{33}(u)\right)\\
  &=&\frac{1}{u-v}\left(\phi(h(-u))-\phi(h(-v))\right) 
\end{eqnarray}
We used the relations (\ref{RLL}) and $[\phi(h(-v)),{L}^{33}(u)]=0$.
For the other relations, the proofs are similar once one remarks, in
particular, that 
\begin{eqnarray}
 {L}^{21}(u){L}^{33}(u)^{-1}&=&{L}^{22}(u){L}^{33}(u)^{-1}
  \phi(f(-u))-\phi(f(-u))-[\phi(f_0),{L}^{22}(u){L}^{33}(u)^{-1}],\label{21}\\
  {L}^{12}(u){L}^{33}(u)^{-1}&=&\phi(e(-u+1))
  -{L}^{22}(u){L}^{33}(u)^{-1}\phi(e(-u+1))
  -[\phi(e_0),{L}^{22}(u){L}^{33}(u)^{-1}],\label{12}\hspace{25pt}\\
 {L}^{31}(u){L}^{33}(u)^{-1}&=&(\phi(f(-u)))^2+\{\phi(f_0),\phi(f(-u))\},
 \label{31}\\
  {L}^{13}(u){L}^{33}(u)^{-1}&=&(\phi(e(-u+1)))^2-\{\phi(e_0),\phi(e(-u+1))\}
 \label{13} .
\end{eqnarray}
The second step consists in proving the surjectivity of $\phi$.
The relations 
(\ref{isoe})-(\ref{isoh}), (\ref{21})-(\ref{13}) and the following 
particular relations, coming from (\ref{C}),
\begin{eqnarray}
C^{22}(u)&=&L^{22}(u-\kappa)L^{22}(u)+L^{32}(u-\kappa)L^{12}(u)
-L^{12}(u-\kappa)L^{32}(u)=\id\\
C^{33}(u)&=&L^{11}(u-\kappa)L^{33}(u)+L^{21}(u-\kappa)L^{23}(u)
-L^{31}(u-\kappa)L^{13}(u)=\id
\end{eqnarray}
constitute a system of nine equations. We can show that these equations are
independent and allow us 
to express all the generators of $\cU(R)$ in terms of $\phi(e_n)$, 
$\phi(h_k)$ and $\phi(f_l)$ ($n,k,l\ge\ZZ_{\ge0}$).
 This proves the surjectivity of $\phi$.
\\
The final step is the proof of the injectivity of $\phi$.
The map $\phi$ preserves the filtration, therefore defines a 
surjective morphism between $gr(\cA^+)$ and $gr(\cU(R))$.
Since the injectivity of the latter morphism is given by the 
proposition \ref{grad}, the injectivity of $\phi$ is proved. 
 \finproof

\medskip\noindent
Note that the RLL relations encode both the commutation relations
\emph{and} the Serre relations.

\medskip\noindent
Let $\overline{\cE}^{\;+}$ and $\overline{\cF}^{\;+}$ be the subalgebras
 of $\cA^+$, without the 
unit, generated by 
$\{e_k, h_l|k,l \ge 0\}$ and $\{f_k, h_l|k,l \ge 0\}$, respectively. 
Let $\cE^+$, $\cH^+$ and $\cF^+$ be the subalgebras of $\cA^+$ 
generated by $e_k$,  $h_k$ and $f_k$ with $k\ge 0$, respectively and 
$\widetilde\cE^+$, $\widetilde\cH^+$ and 
$\widetilde\cF^+$ be the same algebras with the unit.\\ 
\begin{proposition} 
\label{iso_delta}
$\phi$ provides $\cA^+$ with a coalgebra structure given by :
\\
$\bullet$   counit:  
  \begin{eqnarray}
    &&\varepsilon(e(u))=0\;, \qquad   \varepsilon(f(u))=0\;, \qquad
    \varepsilon(h(u))=1\;.\label{counit} 
  \end{eqnarray}
$\bullet$   coproduct :
  \begin{eqnarray}
    \hspace*{-50pt}\Delta(e(u))&\!\!=\!\!& e(u)\otimes\id + \Big(h(u)\otimes
      e(u)+\left[h(u),f_0\right]\otimes({e(u)}^2-\left\{e(u),e_0\right\})
    \Big)\, \Big(\id\otimes \id+X_{12}(u)\Big) \label{coproducte}\\ 
   X_{12}(u) &\!\!=\!\!&\sum_{k>0}(-1)^k\Big(f(u-1)\otimes
      e(u)+({f(u-1)}^2+\left\{f(u-1),f_0\right\}) 
      \otimes({e(u)}^2-\left\{e(u),e_0\right\}) \Big)^k \nonumber\\ 
    \hspace*{-50pt}\Delta(f(u))&\!\!=\!\!& \id\otimes f(u) +
      \Big(f(u)\otimes h(u)+ 
      ({f(u)}^{2}+\left\{f(u),f_0\right\})\otimes
      [h(u),e_0]\Big)\,\Big(\id\otimes \id+Y_{12}(u)\Big)\label{coproductf}\\ 
    Y_{12}(u) &\!\!=\!\!& \sum_{k>0}(-1)^k\Big(f(u)\otimes
      e(u+1)+({f(u)}^{2}+\left\{f(u),f_0\right\}) 
      \otimes({e(u+1)}^{2}-\left\{e(u+1),e_0\right\})
    \Big)^k \nonumber\\   
    \hspace*{-50pt}\Delta(h(u))&\!\!=\!\!& \id \otimes \id + 
\left\{\Delta(e(u)),f_0\otimes
      \id+\id \otimes f_0\right\}\label{coproducth}
  \end{eqnarray}
\end{proposition} 
\textbf{Proof:}
To clarify this proof, we denote $\Delta_A$ (resp. $\Delta_U$) the coproduct
of $\cA^+$ (resp. $U(R)$) and $\eps_A$ (resp. $\eps_U$) the counit
of $\cA^+$ (resp. $U(R)$).\\
We construct $\Delta_A$ thanks to the relation 
$\Delta_U \circ \phi=(\phi \otimes \phi) \circ \Delta_A$. At first, we calculate
$\Delta_U(f(u))$. We begin by
\begin{eqnarray}
  \Delta_U(\phi(f(-u)))&=&\Delta_U\left({L}^{32}(u){L^{33}(u)}^{-1}\right)=
  \Delta_U\left({L}^{32}(u)\right){\Delta_U\left(L^{33}(u)\right)}^{-1} \\
  \mbox{and}\;\Delta_U\left({L}^{33}(u)\right)^{-1}&=&
  \left(\sum_{k=1}^3 ({L}^{3k}(u){L}^{33}(u)^{-1}
    \otimes {L}^{k3}(u){L}^{33}(u)^{-1})({L}^{33}(u)\otimes
    {L}^{33}(u))\right)^{-1}\nonumber \\ 
  &=&({L}^{33}(u)^{-1}\otimes {L}^{33}(u)^{-1})\sum_{n\ge 0}\left(
    \sum_{k=1}^2 {L}^{3k}(u){L}^{33}(u)^{-1}\otimes
  {L}^{k3}(u)L^{33}(u)^{-1}\right)^n\label{delta_33} \hspace{25pt}
\end{eqnarray}
Using the results of the proof of the previous theorem, 
we get 
\begin{eqnarray}
  \Delta_U(\phi(f(-u)))&=&(\phi \otimes \phi)
  \Big[ \id\otimes f(-u) + \Big(f(-u)\otimes h(-u)+
      ({f(-u)}^{2}+\left\{f(-u),f_0\right\}\otimes
      [h(-u),e_0])\Big)\nonumber\\
 &&\hspace{120pt}\times\Big(\id\otimes \id+Y_{12}(-u)\Big)\Big]\nonumber\\
 &=&(\phi \otimes \phi) \circ \Delta_A (f(-u)) \nonumber
\end{eqnarray}
By the injectivity of $\phi$, we find (\ref{coproductf}). 
For $\Delta_A(e(u))$, the proof is similar 
and for $\Delta_A(h(u))$, the equality (\ref{coproducth}) 
is obvious.\\
For the counit, the proof is similar by using $\eps_A=\eps_U\circ\phi$.
\finproof\\
Unfortunately, no explicit formula is  
known for the coproduct in terms of the modes.  \\
Note that (\ref{coproducte})-(\ref{coproducth}) imply that:
  \begin{eqnarray}
    \Delta(e(u))&=& e(u)\otimes\id + h(u)\otimes e(u) 
+ mod(\overline{\cF}^{\;+}\otimes    \cE^+ \cE^+),\\ 
\Delta(f(u)) &=&\id\otimes f(u) + f(u)\otimes h(u) + mod(\cF^+\cF^+\otimes
    \overline{\cE}^{\;+}),\\ 
\Delta(h(u)) &=&h(u)\otimes h(u) 
+ mod(\overline{\cE}^{\;+}\otimes \overline{\cF}^{\;+})
    \label{coproducth2}
  \end{eqnarray}
 To prove (\ref{coproducth2}), we need to calculate the anticommutator of the 
 relation (\ref{coproducth}).\\

\section{The construction of the double $DY(osp(1|2))$}
\setcounter{equation}{0}
In the following, we replace the notations $L(u)$, ${L}^{ij}(u)$, $e(u)$, 
$f(u)$ and $h(u)$ by $L^+(u)$, ${L^+}^{ij}(u)$, $e^+(u)$, $f^+(u)$ 
and $h^+(u)$ respectively. 
\subsection{RTT presentation}
\begin{definition}
  Let $D\cU(R)$ be the associative superalgebra generated by the
  elements $L^{ij}_{(n)}$ 
  $(1\le i,j\le 3$, $n\in\ZZ)$, $\id$ and the defining relations,
  given in terms of formal series  
  $L^+(u)
  =\id_3 +\sum_{i,j=1}^{3} \sum_{n \in \ZZ_{>0}}
  L^{ij}_{(n)} \, u^{-n} \, E_{ij}=\sum_{i,j=1}^{3} {L^+}^{ij}(u) \, E_{ij}$
  and 
  $L^-(u)
  =\sum_{i,j=1}^{3} \sum_{n \in \ZZ_{\le 0}}
  L^{ij}_{(n)} \, u^{-n} \, E_{ij}=\id_3 +\sum_{i,j=1}^{3} {L^-}^{ij}(u) \, E_{ij}$:
  \begin{eqnarray}
    \bullet &&R_{12}\left(u-v\right) \, L^\pm_1(u) \, L^\pm_2(v)
    = L^\pm_2(v) \, L^\pm_1(u) \, R_{12}\left(u-v\right)\label{RL++}\\ 
    \bullet &&R_{12}\left(u-v\right) \, L^+_1(u) \, L^-_2(v) 
    = L^-_2(v) \, L^+_1(u)\, R_{12}\left(u-v\right)\label{RL+-}\\   
    \bullet &&C^\pm(u) = {L^\pm}^t\left(u-\kappa\right) \,
     \, L^\pm(u)=\id_{3}\label{C+-} 
  \end{eqnarray}
  The Hopf algebra structure of $D\cU(R)$ is given by the relations
 (\ref{eq:coproduct})  and (\ref{eq:antipode}) 
 with the substitution $L(u)\,\rightarrow\,L^\pm(u)$.
\end{definition}

\begin{proposition}
The bilinear form $<,>$ between the two subalgebras of $D\cU(R)$,
$\cU^-(R)=\{L^{ij}_{(n)} | n\in\ZZ_{\le0}\}$ with opposite coproduct
and  
$\cU^+(R)=\cU(R)=\{L^{ij}_{(n)} | n\in\ZZ_{>0}\}$ given by:
\begin{eqnarray}
  &&<L^-_1(u), L^+_2(v)>=R_{21}^{-1}\left(v-u\right)
  \mbox{ i.e. }<{L^-}^{ij}(u),{L^+}^{kl}(v)>=
  {\left(R^{-1}\left(v-u\right)\right)}_{ki}^{lj}\label{pair1}\\
  &&<L^-(u),\id_3>=<\id_3,L^+(v)>=<\id_3,\id_3>=\II\label{pair2}
\end{eqnarray}
is a Hopf pairing i.e. satisfies the conditions for $a,b\in\cU^-(R)$
and $\alpha, \beta \in\cU^+(R)$: 
\begin{eqnarray}
  &&<a,\alpha\beta>=(-1)^{[\alpha][\beta]}< \Delta(a), \beta \otimes
  \alpha>\;,\;<ab,\alpha>=<a \otimes b,
  \Delta(\alpha)>\label{pairdelta}\\ 
  &&\varepsilon(a)=<a,\id>,\;\varepsilon(\alpha)=<\id,\alpha>\;,\;
  <S^{-1}(a),\alpha>=<a,S(\alpha)>\label{pairS}\\ 
  &&<a \otimes b,\alpha \otimes \beta>=(-1)^{[b][\alpha]}<a,\alpha><b,\beta>
\end{eqnarray}
\end{proposition}
\textbf{Proof:} The proof for the consistency of (\ref{pair1}) and
(\ref{pair2}) with the defining relations
(\ref{RL++})-(\ref{C+-}) and  
the conditions (\ref{pairdelta}) and (\ref{pairS}) uses the same
methods as in \cite{vlad}. For example : 
\begin{eqnarray}
  &&<L^-_0(w),(R_{12}(u-v)  L^+_1(u) \, L^+_2(v)-L^+_2(v)\, L^+_1(u)
  R_{12}(u-v))> \nonumber\\ 
  &=&R_{12}(u-v)<L^-_0(w),L^+_2(v)><L^-_0(w),L^+_1(u)>
  \nonumber \\
  &&-<L^-_0(w),L^+_1(u)><L^-_0(w),L^+_2(v)>R_{12}(u-v) \\
  &=&R_{12}(u-v)R_{20}^{-1}(v-w)R_{10}^{-1}(u-w)-
  R_{10}^{-1}(u-w)R_{20}^{-1}(v-w)R_{12}(u-v)\hspace{25pt}\\
  &=&0\; \mbox{(due to (\ref{YBE}))}  
\end{eqnarray}
\finproof\\ 
We will show below (see remark after theorem \ref{pairing_B}) 
that this pairing is not degenerate.

\subsection{Drinfel'd realisation}

\begin{definition}
Let $D\cA$ be the associative superalgebra generated 
by the unit $\id$, the even elements h$_k$ ($k\in\ZZ$) and  
the odd elements e$_k$, f$_k$ ($k\in\ZZ$),
gathered in the generating functions 
\begin{eqnarray}
  \label{eq:epemtralala}
  &&e^+(u)=\sum_{k=0}^{\infty} e_k u^{-k-1}\;,
  \qquad
  e^-(u)=-\sum_{k=-\infty}^{-1}e_ku^{-k-1}\;,
  \qquad
  e(u)=e^+(u)-e^-(u)\;,\\
  &&f^+(u)=\sum_{k=0}^{\infty} f_k u^{-k-1}\;,
  \qquad
  f^-(u)=-\sum_{k=-\infty}^{-1}f_ku^{-k-1}\;,
  \qquad
  f(u)=f^+(u)-f^-(u)\;,\\
  &&h^+(u)= \id + \sum_{k=0}^{\infty} h_k u^{-k-1}\;
  \qquad
  h^-(u)= \id - \sum_{k=-\infty}^{-1}h_ku^{-k-1} 
\end{eqnarray}
satisfying the relations
\begin{eqnarray} 
  &&h^\alpha(u)h^\beta(v)=h^\beta(v)h^\alpha(u)\; \mbox{where}\;
  \alpha,\;\beta=\pm,\label{DY_hh}\\ 
  &&\{e(u),f(v)\}=\delta(u-v)\left(h^+(u)-h^-(u)\right)\;
  \mbox{with~}\;\delta(u-v)=\sum_{k=-\infty}^{\infty}u^kv^{-k-1}\\ 
  &&(u-v-1)(2u-2v+1)e(u)h^\pm(v)=(u-v+1)(2u-2v-1)h^\pm(v)e(u),\\
  &&(u-v+1)(2u-2v-1)f(u)h^\pm(v)=(u-v-1)(2u-2v+1)h^\pm(v)f(u),\\
  &&(u-v-1)(2u-2v+1)e(u)e(v)=-(u-v+1)(2u-2v-1)e(v)e(u),\label{Dee}\\  
  &&(u-v+1)(2u-2v-1)f(u)f(v)=-(u-v-1)(2u-2v+1)f(v)f(u)\label{Dff}
\end{eqnarray}
and the supplementary Serre relations  
\begin{eqnarray}  
  && 
  \{e_0,e^\pm(u)\} e^\pm(u) = 
  \frac52 \{e_0^2,e^\pm(u)\} - 2u[e_0^2,e^\pm(u)] + e_0\,
  e^\pm(u) e_0 - [\{e_0,e_1\},e^\pm(u)] \;,
  \label{DY_se}\\   
  && 
  \{f_0,f^\pm(u)\} f^\pm(u) = -
  \frac52 \{f_0^2,f^\pm(u)\} - 2u[f_0^2,f^\pm(u)] - f_0
  f^\pm(u) f_0 - [\{f_0,f_1\},f^\pm(u)] 
  \label{DY_sf} \;.
\end{eqnarray}
The  bialgebra structure of $D\cA$ is given by 
 (\ref{coproducte}), (\ref{coproductf}), (\ref{coproducth}) and (\ref{counit}),
 adding superscripts $\pm$ to $e(u)$, $f(u)$ and $h(u)$. 
\end{definition}
\begin{rmk} $D\cA$ could be alternatively defined by 
the relations (\ref{gen_hh})-(\ref{gen_ff}), 
 adding a superscript $\epsilon $ to the generating functions 
with parameter $u$, and a superscript $\epsilon '$ to the generating functions 
with parameter $v$,
where $\epsilon ,\epsilon'=\pm$.
The Serre relations (\ref{eq:serree}) and (\ref{eq:serref}) with
$\epsilon$ are also valid in $D\cA$, but not sufficient, because they
do not couple enough positive modes with the negative ones. 
\end{rmk} 

The commutation relations in $D\cA$ between $e_k$, $f_k$ and $h_k$
($k\in\ZZ$) are the same as the relations (\ref{h_h})--(\ref{f_f}) 
with, in this case, $k, l \in \ZZ$ and with the following additional 
relations, for $k\in\ZZ$ : 
\begin{eqnarray}
  2[h_{-1},e_k]-2e_{k-1}&=&[h_{-1},e_{k-2}]-\{h_{-1},e_{k-1}\},\label{h-1_e}\\
  2[h_{-1},f_k]-2f_{k-1}&=&[h_{-1},f_{k-2}]+\{h_{-1},f_{k-1}\}\label{h-1_f}
\end{eqnarray}
Similarly to section \ref{sect:A+}, we conjecture that the Serre
relations in terms of modes are of the form (\ref{e_E1})--(\ref{f_F3})
with now $k\in\ZZ$.
As before, this conjecture is supported by explicit computations on
the first orders. 
Moreover, we can define a 
graded algebra for $D\cA$, $grad(D\cA)$, as in proposition
\ref{grad}. Then 
the Serre relations in $grad(D\cA)$ are the relations
(\ref{e_E1})-%
-(\ref{f_F3}), for $k \in \ZZ$,  substituting the
right hand side of the  
equalities by zero. 
Beside this, the expansion of (\ref{DY_se}), (\ref{DY_sf}) in terms of
modes shows that the remaining terms of the Serre relations in $D\cA$ 
have strictly lower degree and are also cubic.  
These results are sufficient to prove the next
theorems.

\null

\noindent
 Let $\cA^+$ and $\cA^-$ be the subalgebras of $D\cA$ generated 
respectively by
 $\{e_n, f_n, h_n | n\in\ZZ_{\ge0}\}$ and 
 $\{e_k, f_k, h_k,E_{-1},F_{-1}| k\in\ZZ_{<0}\}$.
\begin{theorem}
\label{isodouble}
The linear maps
$\Phi:\ D\cA\ \longrightarrow\ D\cU(R)$
and
$\phi^\pm:\ \cA^\pm\ \longrightarrow\ \cU^\pm(R)$
given by
\begin{eqnarray}
 e^\pm(-u)&\longmapsto&{L^\pm}^{33}(u)^{-1}{L^\pm}^{23}(u) \\
 f^\pm(-u)&\longmapsto&{L^\pm}^{32}(u){L^\pm}^{33}(u)^{-1} \\
 h^\pm(-u)&\longmapsto&{L^\pm}^{22}(u){L^\pm}^{33}(u)^{-1}
 +{L^\pm}^{32}(u){L^\pm}^{33}(u)^{-1}{L^\pm}^{23}(u){L^\pm}^{33}(u)^{-1}
\end{eqnarray}
are isomorphisms of bialgebra.
\end{theorem}
\textbf{Proof:} The proof is similar to the one of theorem
\ref{iso} and proposition \ref{iso_delta}. \finproof\\ 
For later convenience, we introduce the following combinations in $D\cA$
(for $k\in \ZZ$):
\begin{eqnarray}
E_{2k+1}&=&\{e_k,e_{k+1}\}-\frac{{e_k}^2}{4} \\
F_{2k+1}&=&\{f_k,f_{k+1}\}-\frac{{f_k}^2}{4}. 
\end{eqnarray}
In particular, $E_{-1}$ and  
$F_{-1}$ will be essential as well as their image in $\cU^-(R)$ :
\begin{eqnarray}
  \Phi:E_{-1}&\longmapsto&-\frac{5}{4}\left
  (\sum_{k\ge 0}(-L^{33}_{(0)})^kL^{23}_{(0)}\right)^2
  -\sum_{k\ge 0}(-L^{33}_{(0)})^kL^{13}_{(0)}\label{E-1}\\
  F_{-1}&\longmapsto&-\frac{5}{4}\left(L^{32}_{(0)}
  \sum_{k\ge 0}(-L^{33}_{(0)})^k\right)^2
  +L^{31}_{(0)}\sum_{k\ge 0}(-L^{33}_{(0)})^k\label{F-1}
\end{eqnarray}
\begin{proposition}
\label{pairingA}
The Hopf pairing $<,>$ between $\cA^-$ and
$\cA^+$ is given by :
\begin{eqnarray}
&&<f^-(u),e^+(v)>=\frac{1}{u-v},\ \ <e^-(u),f^+(v)>
=\frac{-1}{u-v},\label{pair:f-e}\\
&&<h^-(u),h^+(v)>=\frac{(u-v-1 )(2u-2v+1 )}
{(u-v+1 )(2u-2v-1 )}\label{pair:h-h}\\
&&<F_{-1},{e_0}^2>=1,\ \ <E_{-1},{f_0}^2>=1\label{pair:F-e}
\end{eqnarray}
or in terms of generators $(n,k\geq0)$:
\begin{eqnarray}
&&<\id,\id>=1,\;<f_{-k-1},e_n>=- \delta_{n,k},\;
<e_{-k-1},f_n>= \delta_{n,k},\label{pair_mode:f-e}\\
&&<h_{-k-1},h_n>=-\frac{1}{3}\left(2+\left(-\frac{1}{2}\right)^{n-k}\right)
\binom{n}{k}\label{pair_mode:h-h}
\end{eqnarray}
\end{proposition}
\textbf{Proof:} We use the theorem \ref{isodouble} to prove this proposition, 
 for example
\begin{eqnarray}
<e^-(-u),f^+(-v)>&=&<{L^-}^{33}(u)^{-1}{L^-}^{23}(u),
{L^+}^{32}(v){L^+}^{33}(v)^{-1}> \\
&=&<{L^-}^{33}(u)^{-1}\otimes{L^-}^{23}(u),
\Delta({L^+}^{32}(v))\Delta({L^+}^{33}(v))^{-1}> \label{eq3avant}\\
&=&<{L^-}^{23}(u),{L^+}^{32}(v){L^+}^{33}(v)^{-1}>\label{eq3} \\
&=&<{L^-}^{22}(u),{L^+}^{33}(v)>^{-1}
<{L^-}^{23}(u),{L^+}^{32}(v)>=\frac{1}{u-v} 
\end{eqnarray}
The difficult point lies in the step between equalities 
(\ref{eq3avant}) and (\ref{eq3}). 
It is done using the explicit form (\ref{delta_33}) 
and showing that only the first term of the sum contributes 
to the pairing.\\
For $<f^-(u),e^+(v)>$, the  proof is similar.\\
The identity $<h^-(u),h^+(v)>=<h^-(u),\{e^+(v),f_0\}+\id>$
 and the two previous results allow us to obtain the relation
 (\ref{pair:h-h}).\\
The pairings $<F_{-1},{e_0}^2>$ and $<E_{-1},{f_0}^2>$ are calculated 
thanks to the explicit expressions (\ref{E-1}) and
 (\ref{F-1}).\\
To find the explicit form (\ref{pair_mode:h-h}), we remark that 
($n\geq0$):
\begin{eqnarray}
&&<h^-(u),h_0>=1,\qquad <h^-(u),h_1>=u+\frac{1}{2} \\
&&<h^-(u),h_{n+2}>-\left(2u+\frac{1}{2}\right)<h^-(u),h_{n+1}>
+(u+1)\left(u-\frac{1}{2}\right)<h^-(u),h_n>=0.\ 
\end{eqnarray}
A trivial induction shows that
 $<h^-(u),h_{n+2}>=\frac{1}{3}\left(u-\frac{1}{2}\right)^{n+2}
+\frac{2}{3}\left(u+\frac{1}{2}\right)\left(u+1\right)^{n+1}
+\frac{1}{3}\left(u+1\right)^{n+1}$ which gives the result.
\finproof\\

\subsection{Dual bases}

Now, we look for bases of $\cA^+$ and $\cA^-$.
 Let $\overline{\cE}^{\;-}$ and $\overline{\cF}^{\;-}$ be the subalgebras
 of $\cA^-$, generated by 
$\{e_k, E_{-1}, h_l|k,l < 0\}$ and $\{f_k, F_{-1}, h_l|k,l < 0\}$ respectively. 
Let $\cE^-$, $\cH^-$ and $\cF^-$ be the subalgebras of $\cA^-$, 
without the unit, generated by $\{e_k,E_{-1}|k<0\}$,  $\{h_k|k<0\}$ 
and $\{f_k,F_{-1}|k<0\}$, respectively and $\widetilde\cE^-$, 
$\widetilde\cH^-$ and $\widetilde\cF^-$ be the same algebras 
with the unit.
\begin{proposition}
\label{decompo_trian}
$(i)$ Let $a_\pm \in \cA^\pm$, then $a_+ \in \widetilde\cE^+
\widetilde\cH^+\widetilde\cF^+$ 
and $a_- \in \widetilde\cF^-\widetilde\cH^-\widetilde\cE^-$.\\
$(ii)$ $\forall e_\pm \in \widetilde\cE^\pm,
 h_\pm \in \widetilde\cH^\pm, f_\pm \in \widetilde\cF^\pm,$
\begin{equation} 
<f_-h_-e_-,e_+h_+f_+>
=(-1)^{[e_+][e_-]}<f_-,e_+><h_-,h_+><e_-,f_+>
\end{equation}
\end{proposition}
\textbf{Proof:} $(i)$ We use a proof similar to the one given in \cite{Rosso}.
We first consider $\cA^+$.
As $\{e_k, h_l, f_m |k,l,m\ge 0\}$ is a generating set of $\cA^+$, 
it is enough to prove that any monomial
$\prod_{i=0}^j {e_{k_i}}^{\alpha_i}{h_{l_i}}^{\beta_i}{f_{m_i}}^{\gamma_i}$,
with $j, k_i,\alpha_i, l_i,\beta_i, m_i,\gamma_i \in \ZZ_{\ge 0}$, 
is a linear combination of elements belonging to 
$\widetilde\cE^+\widetilde\cH^+\widetilde\cF^+$. 
We make an induction on the degree 
$p=\sum_{i=0}^j(\alpha_i+\beta_i+\gamma_i)$ of 
such monomial.\\ 
For $p=1$, the assertion is obvious.\\
Let assume the assertion is true for $p\ge1$ and consider an element 
$\prod_{i=0}^j {e_{k_i}}^{\alpha_i}{h_{l_i}}^{\beta_i}
{f_{m_i}}^{\gamma_i}$ 
such that $\sum_{i=0}^j(\alpha_i+\beta_i+\gamma_i)=p+1$. 
The last $p$ generators can be ordered using the induction
hypothesis on $p$. Then, three
 cases are possible depending on the first 
element:
\begin{itemize}
\item It belongs to $\widetilde\cE^+$. Then, the element is ordered.
\item It belongs to $\widetilde\cH^+$. If the second generator 
belongs to $\widetilde\cH^+$ or $\widetilde\cF^+$, the assertion for $p+1$ 
is proven. It remains the case where the second generator (say $e_k$)
belongs to $\cE^+$. We make 
an induction on the index, $l$, of the 
first generator $h_l$.
Let $a_{p-1}$ be the ordered product of
the last $p-1$ generators.
\\
If $l=0$ then 
\begin{equation}
h_0e_ka_{p-1}=e_kh_0a_{p-1}+e_ka_{p-1}\;
\mbox{due to}\; (\ref{h0_e})
\end{equation}
As the induction on $p$ allows us to order $h_0a_{p-1}$, 
we can order this element.\\
If $l=1$ then 
\begin{equation}
h_1e_ka_{p-1}=e_kh_1a_{p-1}+e_{k+1}a_{p-1}
+h_0e_ka_{p-1}+e_kh_0a_{p-1}\;\mbox{due to}\;
 (\ref{h1_e})
 \end{equation}
Using the case $l=0$ for $h_0e_ka_{p-1}$ and the induction hypothesis on 
$p$ for $h_1a_{p-1}$ and $h_0a_{p-1}$, we can order 
this element.\\
Let $l\ge2$ and assume that for $l-2$ and $l-1$ the elements 
can be ordered. 
Then, using the commutation relations (\ref{h_e}), one gets 
\begin{equation}
h_le_ka_{p-1}=
(-h_{l-2}e_{k+2}+h_{l-1}e_k/2-h_{l-2}e_{k+1}/2
+h_{l-2}e_k/2+2h_{l-1}e_{k+1})a_{p-1}+b_{p+1}
\end{equation}
where $b_{p+1}$ can be ordered thanks to the induction on $p$. 
The other elements are ordered by the induction on 
$l-2$ and $l-1$.\\
This ends the induction on $l$.
\item It belongs to $\widetilde\cF^+$. We denote it $f_m$. 
If the second element belongs to 
$\widetilde\cF^+$, 
the assertion for $p+1$ is proven. 
If the second element belongs to $\widetilde\cH^+$, we prove 
the assertion for $p+1$ analogously to
the previous case, using (\ref{h0_f}), (\ref{h1_f}) and (\ref{h_f}).
Finally, if the second element $e_k$ belongs to $\widetilde\cE^+$, then 
\begin{equation}
f_me_ka_{p-1}=-e_kf_ma_{p-1}+h_{m+k}a_{p-1}\;
\mbox{due to}\;(\ref{e_f})
\end{equation}
and $f_ma_{p-1}$, $h_{m+k}a_{p-1}$ are ordered thanks 
to the hypothesis on $p$.
\end{itemize}
This ends the induction on $p$ and $(i)$ 
 is proven for $\cA^+$.\\
For $\cA^-$, the proof is almost similar. However, an additional
difficulty appears because of  
exchange relations between $e_{-k}$ and $h_{-1}$.
The relation (\ref{h-1_e}) allows us to order $\hat e_{-l}
\equiv 2e_{-l}-e_{-l-1}-e_{-l-2}$ 
($l\ge1$) and $h_{-1}$. Fortunately, $e_{-k}$ can be expressed 
in terms of 
$\hat e_{-l}$ :
\begin{equation}
\forall k\ge 1,\;e_{-k}=\frac{1}{3}\sum_{l=k}^{+\infty}
\left(1-{\left(-\frac{1}{2}\right)}^{l-k+1}\right)\hat e_{-l}
\end{equation}
Therefore, we can order $e_{-k}$ and $h_{-1}$.
Likewise, $e_{-k}$ and $h_{-n}$ can be ordered. 
\\
$(ii)$ Let $\bar{f}_- \in \widetilde{\cF}^{-}\widetilde{\cH}^{-}, 
\bar{e}_+\in \widetilde{\cE}^{+}\widetilde{\cH}^{+},
 e_- \in \cE^-$ and $f_+ \in \cF^+$. One computes
\begin{eqnarray}
<\bar{f}_-\;e_-,\bar{e}_+\;f_+>&=&<\bar{f}_-\otimes
 e_-,\Delta(\bar{e}_+)\Delta(f_+)> \\
&=&<\bar{f}_-\otimes e_-,(\bar{e}_+\otimes \id
+mod(\cA^+\otimes \overline{\cE}^{\;+}))(\id \otimes f_+ 
+mod( \cF^+\otimes \cA^+))> \hspace{25pt}\\
&=&(-1)^{[\bar{e}_+][e_-]}<\bar{f}_-,\bar{e}_+><e_-, f_+> 
\end{eqnarray}
using the following identities:
\begin{eqnarray}
<e_-, \overline{\cE}^{\;+} f_+> 
&=&(-1)^{[f_+][\overline{\cE}^{\;+}]}
<\Delta(e_-), f_+\otimes\overline{\cE}^{\;+}>\nonumber\\
&=&(-1)^{[f_+][\overline{\cE}^{\;+}]}<e_-\otimes \id+mod(\cA^-\otimes \cE^-),
f_+\otimes\overline{\cE}^{\;+}>=0 \\
<\bar{f}_-,\bar{e}_+\cF^+> 
&=&(-1)^{[\bar{e}_+][\cF^+]}
<\Delta(\bar{f}_-),\cF^+\otimes \bar{e}_+>\nonumber\\
&=&(-1)^{[\bar{e}_+][\cF^+]}
 <\id\otimes\bar{f}_-+mod(\overline{\cF}^{\;-}\otimes \cA^-),
\cF^+\otimes \bar{e}_+>=0 \\
<\bar{f}_-,\cA^+\cF^+> 
&=&(-1)^{[\cF^+][\cA^+]} 
<\Delta(\bar{f}_-),\cF^+\otimes\cA^+>\nonumber\\
&=&(-1)^{[\cF^+][\cA^+]} 
<\id\otimes\bar{f}_-+mod(\overline{\cF}^{\;-}\otimes \cA^-),
\cF^+\otimes\cA^+>=0  \hspace{35pt}
\end{eqnarray}
Let $f_-\in \cF^-$, $h_-\in \cH^-$, 
$e_+\in \cE^+$ and $h_+\in \cH^+$ such that $\bar{f}_-=f_-h_-$
and $\bar{e}_+=e_+h_+$.
 Then, we prove analogously
that $<\bar{f}_-,\bar{e}_+>=<f_-h_-,e_+h_+>=<f_-,e_+><h_-,h_+>$ using 
$\Delta(h_+)=\id\otimes h_+ 
+\sum_{i}h_i\otimes h'_i
+mod(\cE^+{\cH}^+\otimes{\cH}^+{\cF}^+)$ 
where $h_i$ and $h'_i \in \cH^+$.
Remarking  that $[\bar{e}_+]=[e_+]$ and that for the unit the theorem is obvious,
 we prove the second assertion
 of the theorem.
\finproof\\
\begin{rmk} The point $(ii)$ of the previous theorem shows that the dual
of $\widetilde\cE^-$ (resp. $\widetilde\cH^-$, $\widetilde\cF^-$) 
is $\widetilde\cF^+$ (resp. $\widetilde\cH^+$, $\widetilde\cE^+$ ).
In particular, one has $<\widetilde\cE^-,\widetilde\cE^+>=0=
<\widetilde\cE^-,\widetilde\cH^+>$, and the same relations changing 
$\widetilde\cE$ by $\widetilde\cF$.
\end{rmk}
\begin{theorem}
Bases of $\widetilde\cE^+$, $\widetilde\cF^-$, $\widetilde\cF^+$ 
and $\widetilde\cE^-$ are respectively 
\begin{eqnarray}
  \cB_{E^+}&=&\{{e_0}^{a_0}{E_1}^{b_0}{e_1}^{a_1}{E_3}^{b_1}
  \ldots
  {e_k}^{a_k}{E_{2k+1}}^{b_k} \ldots
  |a_0, a_1\ldots, b_0, b_1\ldots\in
  \ZZ_{\ge 0}\}, \\ 
  \cB_{F^-}&=&\{{F_{-1}}^{b_0}{f_{-1}}^{a_0}{F_{-3}}^{b_1}{f_{-2}}^{a_1}
  \ldots {F_{-2k-1}}^{b_k}{f_{-k-1}}^{a_k} \ldots
  |a_0, a_1\ldots, b_0, b_1\ldots\in \ZZ_{\ge 0}\}, \\
  \cB_{F^+}&=&\{{F_{2k+1}}^{b_k}{f_k}^{a_k}
  \ldots{F_3}^{b_1}{f_1}^{a_1}{F_1}^{b_0}{f_0}^{a_0} \ldots
  |a_0, a_1\ldots, b_0, b_1\ldots\in \ZZ_{\ge 0}\}, \\
  \cB_{E^-}&=&\{{e_{-k-1}}^{a_k}{E_{-2k-1}}^{b_k}
  \ldots{e_{-2}}^{a_1}{E_{-3}}^{b_1}{e_{-1}}^{a_0}{E_{-1}}^{b_0} \ldots
  |a_0, a_1\ldots, b_0, b_1\ldots\in \ZZ_{\ge 0}\}.
\end{eqnarray}
\end{theorem} 
\textbf{Proof:} To prove that $\cB_{E^+}$, $\cB_{F^-}$, 
$ \cB_{F^+}$ and $\cB_{E^-}$ 
generate $\widetilde\cE^+$, $\widetilde\cF^-$, $\widetilde\cF^+$ 
and $\widetilde\cE^-$, respectively,
the methods are the same as the ones used in the proof 
of the assertion $(i)$ of proposition \ref{decompo_trian}, 
using the commutation relations between $e_k$ and $e_l$ 
as well as the proved results on the Serre relations.\\
The independence of the set of generators is given by their 
 independence in the corresponding graded algebras which is 
 obvious.
\finproof\\

\begin{theorem}
\label{pairing_B}
For any element $b_+$ of $\cB_{E^+}$ (resp. $\cB_{F^+}$), there is one
and only one element $b_-$ in $\cB_{F^-}$ (resp. $\cB_{E^-}$) such that
the pairing $<b_-,b_+>$ does not vanish. They are given by:  
\begin{eqnarray}
  &&<\prod_{n=0}^{\stackrel{k}{\longrightarrow}}
   {F_{-2n-1}}^{b_n}{f_{-n-1}}^{2a_n+c_n} 
    ,\prod_{m=0}^{\stackrel{k}{\longrightarrow}}
    {e_m}^{2b_m+c_m}{E_{2m+1}}^{a_m}>
  =(-1)^{\sum_{0\le n < m \le k} c_n c_m}
  \prod_{l=0}^k (-1)^{c_l}a_l!\; b_l! \hspace{25pt} \label{pairBE}\\
   &&<\prod_{n=0}^{\stackrel{k}{\longleftarrow}}
  {e_{-n-1}}^{2a_n+c_n}{E_{-2n-1}}^{b_n},
  \prod_{m=0}^{\stackrel{k}{\longleftarrow}}
  {F_{2m+1}}^{a_m}{f_m}^{2b_m+c_m}>
  =(-1)^{\sum_{0\le n < m \le k} c_n c_m}
  \prod_{l=0}^k a_l!\;b_l!\label{pairBF}\\
&&\mbox{where}\;k \in \ZZ_{\ge 0},\; a_k,b_k\in \ZZ_{\ge 0},\;
c_k\in\{0,1\},\;\prod_{n=0}^{\stackrel{k}{\longrightarrow}}
e_n=e_0\,e_1\ldots\, e_{k-1}\,e_k\; 
\mbox{and}\;
\prod_{n=0}^{\stackrel{k}{\longleftarrow}}e_n=e_k\,e_{k-1}
\ldots\, e_1\,e_0. 
\nonumber
\end{eqnarray}
\end{theorem} 
\textbf{Proof:} 
We first consider the pairing $<a_-,e_k>$ for $a_-\in\cB_{F^-}$.
If $a_-$ has degree at least 2, 
there are $a'_-,a''_- \in \cB_{F^-}$ such that $a_-=a'_-a''_-$.
Using (\ref{pair_mode:f-e}), we get
\begin{eqnarray*}
  <a_-,e_k> &=& <a'_-\otimes a''_-,\Delta(e_k)>
  \\
  &=&<a'_-\otimes a''_-,e_k\otimes \id+\id\otimes e_k
  +\sum_{l=0}^{k-1}h_l\otimes e_{k-l-1}+ 
  mod(\overline{\cF}^{\;+}\otimes \cE^+\cE^+)>\\
  &=&
  <a'_-\otimes a''_-,e_k\otimes \id+\id\otimes e_k
  +\sum_{l=0}^{k-1}h_l\otimes e_{k-l-1}>=0,
\end{eqnarray*}
This shows that the pairing of $<a_-,e_k>$ is different from zero only for
$a_-=f_{-k-1}$ and given by (\ref{pair_mode:f-e}).\\
Consider now  $<a_-,e_k^2>$. If $a_-=f_{-m}$ has degree 1, one computes:
\begin{equation}
<f_{-m},e_k^2>=<\Delta(f_{-m}),e_k\otimes e_k>=<f_{-m}\otimes \id+\id\otimes
f_{-m}
+\sum_{l=0}^{k-1}h_l\otimes f_{-m-l-1},e_k\otimes e_k>=0
\end{equation}
For $a_-=a'_-a''_-$ of degree at least 2, one has
\begin{eqnarray}
<a_-,{e_k}^2> &=& <a'_-\otimes a''_-,{\Delta(e_k)}^2>
=<a'_-\otimes a''_-,\sum_{l=0}^{k-1}[e_k,h_l]\otimes e_{k-l-1}> \\
&=&<a'_-\otimes a''_-,-e_k\otimes e_{k-1}
+\sum_{l=1}^{k-1}\Psi_l(e_k,e_{k+1},\ldots,e_{k+l})\otimes e_{-l-1}> 
\end{eqnarray}
where $\Psi_l(e_k,e_{k+1},\ldots,e_{k+l})$ is a linear combination 
of $e_k,e_{k+1},\ldots,e_{k+l}$.
Then, for $k\ge 1, <a_-,{e_k}^2>$ is equal $1$ for $a'_-=f_{-k-1}$ 
and $a''_-=f_{-k}$ (i.e. for $a_-=F_{-2k-1}$) and $0$ otherwise. 
For $k=0$, by (\ref{pair:F-e}) and the previous calculation,
 we know that the pairing of ${{e_0}^2}$ does not vanish only with $F_{-1}$.\\
Similarly, $<{f_{-k-1}}^2,a_+>$ is equal to 1 if $a_+=E_{2k+1}$, and to
0 in the other cases.\\
Now, we show by induction that $<{F_{-2k-1}}^{b},{e_k}^{2b}>=b!$ 
and $<{F_{-2k-1}}^{b}f_{-k-1},{e_k}^{2b+1}>=-b!$ and 
that the other pairings with ${e_k}^{2b}$ or ${e_k}^{2b+1}$ are zero. 
We assume these assertions for $b<b_0$.
\begin{eqnarray}
&&<a_-,{e_k}^{2b_0}>=<a'_-\otimes a''_-,\left({e_k}^2\otimes\id
+\id\otimes {e_k}^2
+\sum_{l=0}^{k-1}\Psi_l(e_k,e_{k+1},\ldots,e_{k+l})
\otimes e_{k-l-1}\right)^{b_0}>\hspace{25pt} \nonumber\\
&=&<a'_-\otimes a''_-,\sum_{p=0}^{b_0}\binom{b_0}{p}{e_k}^{2p}
\otimes {e_k}^{2(b_0-p)}>=
<a'_-\otimes a''_-,\sum_{p=1}^{b_0-1}\binom{b_0}{p}{e_k}^{2p}
\otimes {e_k}^{2(b_0-p)}> 
\end{eqnarray}
Due to the hypothesis, this pairing is non zero only if
$a'={F_{-2k-1}}^{p}$ and  $a''={F_{-2k-1}}^{b_0-p}$ and, in this case,
is equal to  
\begin{equation}
\binom{b_0}{p}<{F_{-2k-1}}^{p},{e_k}^{2p}><{F_{-2k-1}}^{b_0-p},
{e_k}^{2(b_0-p)}>=\binom{b_0}{p}p!\;(b_0-p)!=b_0!
\end{equation} 
A similar result is proven for ${e_k}^{2b_0+1}$ and then by induction on
$b$, the assertions are proven.\\ 
Similarly, $<{f_{-k-1}}^{2a},{E_{2k+1}}^a>=a!$ and
$<{f_{-k-1}}^{2a+1},e_k{E_{2k+1}}^a>=-a!$.\\ 
We can sum up all these results by
\begin{eqnarray}
  \forall
  a,b\in\ZZ_{\ge0},c\in\{0,1\},
  &&<{F_{-2k-1}}^b{f_{-k-1}}^{2a+c},{e_k}^{2b+c}{E_{2k+1}}^a>
  =(-1)^c\,a!\;b!
\end{eqnarray}
and all other pairings with ${e_k}^{2b+c}{E_{2k+1}}^a$ are zero. 
Similarly, we show that all other pairings with
${F_{-2k-1}}^b{f_{-k-1}}^{2a+c}$ are also zero. 
\begin{eqnarray*}
  &&<\prod_{n=0}^{\stackrel{l}{\longrightarrow}}
  {F_{-2n-1}}^{b'_n}{f_{-n-1}}^{2a'_n+c'_n}  
  ,\prod_{m=0}^{\stackrel{k}{\longrightarrow}}
  {e_m}^{2b_m+c_m}{E_{2m+1}}^{a_m}> \\
  &=&<{F_{-1}}^{b'_0}{f_{-1}}^{2a'_0+c_0}\otimes
  \prod_{n=1}^{\stackrel{l}{\longrightarrow}}
  {F_{-2l-1}}^{b'_l}{f_{-l-1}}^{2a'_l+c'_l}  
  ,\Delta({e_0}^{2b_0+c_0}{E_1}^{a_0})
  \Delta(\prod_{m=1}^{\stackrel{k}{\longrightarrow}}
  {e_m}^{2b_m+c_m}{E_{2m+1}}^{a_m})>
    \\ 
  &=&(-1)^{c_0(c'_1+\dots+c'_l}\delta_{a'_0,a_0}\delta_{b'_0,b_0}
  \delta_{c'_0,c_0}(-1)^{c_0}a_0!\,b_0!
  <\prod_{n=1}^{\stackrel{l}{\longrightarrow}}
  {F_{-2l-1}}^{b'_l}{f_{-l-1}}^{2a'_l+c'_l} 
  ,\prod_{m=1}^{\stackrel{k}{\longrightarrow}} 
  {e_m}^{2b_m+c_m}{E_{2m+1}}^{a_m}>    \hspace{35pt}
\end{eqnarray*}
Repeating this calculus $k$ times, we prove (\ref{pairBE}) of the theorem.\\
(\ref{pairBF}) is proven along the same lines.
\finproof\\
\begin{rmk} Since $\widetilde\cH^+$ is Abelian, one of its basis 
 is $\{{h_0}^{a_0}{h_1}^{a_1}\ldots
 {h_k}^{a_k}\ldots|a_0,a_1,\ldots\in  \ZZ_{\ge 0}\}$.
  In addition, the pairing restricted of the subalgebras 
  $\widetilde\cH^-$ and $\widetilde\cH^+$
  is not degenerated. 
\end{rmk}
\begin{rmk} A corollary of the previous results is that the pairing
between $\cA^-$ and $\cA^+$ is not degenerated. Then,
thanks to the isomorphisms $\phi^\pm$ and $\Phi$, neither is
the pairing
between $\cU^-(R)$ and $\cU^+(R)$.
\end{rmk}

\begin{theorem}
\label{doubleR}
$\cU^+(R) \otimes \cU^-(R)$ is the quantum double of $\cU^+(R)$ with
 the multiplication between $\cU^+(R)$ and $\cU^-(R)$ defined by
 (\ref{RL+-}). Thus, it
 is isomorphic, as a Hopf algebra, to the quantum double of
 $Y(osp(1|2))$, denoted $DY(osp(1|2))$. \\
 Similarly, $\cA^+(R) \otimes \cA^-(R)$  is the quantum double of $\cA^+(R)$.
\end{theorem}
\textbf{Proof:} From $(\Delta\otimes \id)\Delta(L^\pm(u))
=L^\pm(u)\dot\otimes L^\pm(u)\dot\otimes L^\pm(u)$, 
the cross-multiplication in a quantum double is defined by  
\begin{eqnarray}
  L_2^-(v)L_1^+(u) &=&
  <S(L_2^-(v)),L_1^-(u)> L_1^+(u) L_2^-(v) <L_2^-(v),L_1^+(u)> \\ 
  &=&R_{12}(u-v)L_1^+(u)L_2^-(v)R_{12}^{-1}(u-v) 
\end{eqnarray}
which is equivalent to (\ref{RL+-}).\\
The other assertions are obvious.
\finproof\\
 
\section{Universal R-matrix}
\setcounter{equation}{0}

\subsection{Construction of the universal R-matrix}
We express the universal R-matrix of double 
super Yangian $DY(osp(1|2))$ according to the generators of 
Drinfel'd basis. 
Since $D\cA$ is the quantum double of $\cA^+$, it 
admits a canonical universal R-matrix given by
$\cR=\sum x_i \otimes x^i$ where $\{x_i,i\in\NN\}$
 is the basis of $\cA^+$ and $\{x^i,i\in\NN\}\in \cA^-$ is the 
dual basis
  (i.e. $<x^i,x_j>=\delta^i_j$). 
 Therefore, thanks to the explicit expression of the pairing, we have 
 the following result.

\begin{theorem} 
The universal R-matrix can be factorised as
 \begin{equation}
 \cR=\cR_E\cR_H\cR_F
 \end{equation}
 where $\cR_E\in \widetilde\cE^+ \otimes  \widetilde\cF^-$, 
 $\cR_H \in  \widetilde\cH^+ \otimes \widetilde\cH^-$ 
 and $\cR_F \in \widetilde\cF^+\otimes\widetilde\cE^- $.\\
  The explicit expressions of the universal factors $\cR_E$ and $\cR_F$ are
  \begin{equation}
  \cR_E=\prod_{i\ge0}^{\longrightarrow}\left[
    \exp\left({e_i}^2\otimes F_{-2i-1}\right) (\id\otimes\id
     -e_i\otimes f_{-i-1}) 
    \exp\left(E_{2i+1}\otimes {f_{-i-1}}^2\right)\right]
    \end{equation}
  \begin{equation}\cR_F=\prod_{i\ge0}^{\longleftarrow}\left[
    \exp\left(F_{2i+1}\otimes {e_{-i-1}}^2\right)
    (\id\otimes\id+f_i\otimes e_{-i-1})
  \exp\left({f_i}^2\otimes E_{-2i-1}\right) 
\right]
\end{equation}
\end{theorem} 
 \textbf{Proof:} The factorisation of the the universal R-matrix 
 is involved by the relation $(ii)$
of the proposition (\ref{decompo_trian}). 
In addition, to prove the expression of $\cR_{E}$,
we expand the exponentials and the products
\begin{eqnarray}
\cR_{E}&=&\prod_{i\ge0}^{\longrightarrow}
\sum_{{a_i},{b_i}\ge 0}\frac{1}{a_i\,!{b_i}\,!}\left({e_i}^{2a_i}
\otimes {F_{-2i-1}}^{a_i}\right)
 (\id\otimes\id-e_i\otimes f_{-i-1})
\left({E_{2i+1}}^{b_i}\otimes {f_{-i-1}}^{2{b_i}}\right)\\
&=&\prod_{i\ge0}^{\longrightarrow}
 \sum_{
 {a_i},{b_i}\ge 0,\atop
 c_i\in\{0,1\}}\frac{(-1)^{c_i}}{a_i\,!{b_i}\,!}
 \left({e_i}^{2a_i+c_i}{E_{2i+1}}^{b_i}
 \otimes {F_{-2i-1}}^{a_i}{f_{-i-1}}^{2{b_i}+c_i}\right)\\
&=&\left(
 \sum_{
 {a_0},{b_0}\ge 0,\atop
 c_0\in\{0,1\}}
 \sum_{
 {a_1},{b_1}\ge 0,\atop
 c_1\in\{0,1\}}\ldots
 \right)
 \frac{(-1)^{c_0+c_1+\cdots}}
 {a_0\,!b_0\,!a_1\,!b_1\,!\ldots}(-1)^{c_0c_1+\cdots}\nonumber\\
 &&\times\left({e_0}^{2a_0+c_0}{E_{1}}^{b_0}
 {e_1}^{2a_1+c_1}{E_{3}}^{b_1}\ldots
 \otimes {F_{-1}}^{a_0}{f_{-1}}^{2{b_0}+c_0}
 {F_{-3}}^{a_1}{f_{-2}}^{2{b_1}+c_1}\ldots\right) 
\end{eqnarray}
Therefore, $\cR_{E}$ can be written $\sum_{x_i\in\cB_{E^+}}x_i \otimes x^i$
and $<x^i,x_j>=\delta^i_j$ due to (\ref{pairBE}).
The proof to find the explicit form of $\cR_{F}$ is similar.
\finproof\\
\begin{theorem}\label{theor:RH}
The factor $\cR_{H}$ of the universal R-matrix of $DY(osp(1|2))$ is given
by
\begin{equation}
\cR_{H} = \prod_{n \ge 0} \exp \left\{ \sum_{i \ge 0} \Big( \frac{d}{du} \,
K_{-}(u) \Big)_{i} \otimes \Big( C(T^{1/2}) K_{+}(v+3n+\three) \Big)_{-i-1}
\right\}
\label{eq:RH}
\end{equation}
where $K_{\pm}(u) = \ln h^{\pm}(u)$, $C(q) = q+1+q^{-1}$, $T$ is the
shift operator: $T f(u) = f(u+1)$ and $(\psi(u))_{i} = \psi_{i}$ for any
function $\psi(u) = \sum_{i} \psi_{i} u^{-i-1}$.
\end{theorem}
\textbf{Proof:} The proof is inspired by the results exposed in \cite{KT}.
Starting from the pairing (\ref{pair:h-h}), a direct calculation shows
that
\begin{equation}
< K_{-}(u),K_{+}(v) > \;=\; \ln \,
\frac{(u-v-1)(2u-2v+1)}{(u-v+1)(2u-2v-1)}
\end{equation}
{from} which it follows
\begin{eqnarray}
< \frac{d}{du} \, K_{-}(u),K_{+}(v) > &=& \frac{1}{u-v-1} +
\frac{1}{u-v+1/2} - \frac{1}{u-v+1} - \frac{1}{u-v-1/2} \\
&=& (T^{-1} + T^{1/2} - T - T^{-1/2}) \; \frac{1}{u-v}
\end{eqnarray}
Therefore one obtains
\begin{equation}
< \frac{d}{du} \, K_{-}(u) , (T^{-1} + T^{1/2} - T - T^{-1/2})^{-1}
K_{+}(v) > \;=\; \frac{1}{u-v}
\end{equation}
The formal inversion of the operator $(T^{-1} + T^{1/2} - T -
T^{-1/2})^{-1}$ is given by
\begin{equation}
(T^{-1} + T^{1/2} - T - T^{-1/2})^{-1} = \sum_{n \ge 0} T^{3n+2} +
T^{3n+3/2} + T^{3n+1}
\label{eq:invT}
\end{equation}
Let $B(q)$ be the $q$-analog of the symmetrised Cartan matrix 
of $osp(1|2)$. We define $C(q)$ by the relation\footnote{The 
presence of $q^{1/2}$ instead of
$q$ in the definition of $C(q)$ is due to the normalisation of the
fermionic simple root.} 
$\displaystyle B(q)^{-1} = \frac{1}{[2\kappa]_{q^{1/2}}} \; C(q)$, $C(q)$
is a matrix with polynomial entries in $q$ and $q^{-1}$ and positive
coefficients. One gets $C(q) = q+1+q^{-1}$. It follows that
\begin{equation}
\sum_{n \ge 0} < \frac{d}{du} \, K_{-}(u) , C(T^{1/2}) \; K_{+}(v+3n+3/2) >
\;=\; \frac{1}{u-v}
\label{eq:invT2}
\end{equation}
Since the pairing (\ref{eq:invT2}) exhibits a duality relation in diagonal
form, one gets immediately the expression (\ref{eq:RH}) for the universal
factor $\cR_{H}$.
\finproof\\

\subsection{Evaluated $R$-matrix}

\begin{proposition}
Let $\pi$ be the fundamental 3-dimensional representation of $osp(1|2)$
with representation space $V$ and $V_{z}$ a $\CC$-module. Then $\pi_{z}$
such that
\begin{eqnarray}
  \pi_{z} : D\cA & \longrightarrow & V_{z} \otimes V \nonumber \\
  e_n &\longmapsto& z^nE_{12}+{z'}^nE_{23} \\
  f_n &\longmapsto&  z^nE_{21}-{z'}^nE_{32} \\
  h_n &\longmapsto& z^nE_{11}+(z^n-{z'}^n)E_{22}-{z'}^nE_{33}
\end{eqnarray}
is an evaluation representation of the double super Yangian $DY(osp(1|2))$
for $z'=z+\half$.
\end{proposition}
\textbf{Proof:} The image by $\pi_z$ of the elements of $D\cA$ have to
satisfy the commutation relations 
 (\ref{DY_hh})-(\ref{DY_sf}).
For example, we prove for (\ref{Dee}). We use that $\pi_{z}(e(u))
=\delta(z-u)E_{12}+\delta(z'-u)E_{23}$.
Then, one gets :
\begin{eqnarray}
&& (u-v-1)(2u-2v+1)\pi_{z}(e(u))\pi_{z}(e(v))
+(u-v+1)(2u-2v-1)\pi_{z}(e(v))\pi_{z}(e(u))\nonumber \\
&=&[(u-v-1)(2u-2v+1)\delta(z-u)\delta(z'-v)
+(u-v+1)(2u-2v-1)\delta(z-v)\delta(z'-u)] {E_{13}}\nonumber\\
&=&[(z-z'-1)(2z-2z'+1)+(z'-z+1)(2z'-2z-1) {E_{13}}=0
\end{eqnarray}
The other commutation relations are proven analogously.
\finproof\\
\begin{theorem}
Let $\pi_{z}$ and $\pi_{w}$ be two fundamental evaluation
representations, then
\begin{equation}
(\pi_z\otimes\pi_w)\cR=R_{12}(z-w)
\end{equation}
where the $R$-matrix $R_{12}(z)$ is given by (\ref{matrixR}).
\end{theorem}
\textbf{Proof:}
\begin{eqnarray*}
&& \hspace*{-10mm} (\pi_z \otimes \pi_w)\cR_E \;=\; \id_3 \otimes \id_3
+\sum_{j>i\ge 0}\pi_z(e_ie_j) \otimes \pi_w(f_{-i-1}f_{-j-1}) \nonumber \\
&& \hspace*{-5mm} + \; \sum_{i\ge 0}\left[\pi_z({e_i}^2)\otimes
\pi_w(F_{-2i-1}) -\pi_z(e_i) \otimes \pi_w(f_{-i-1}) +\pi_z(E_{2i+1})
\otimes \pi_w({f_{-i-1}}^2)\right] \nonumber \\
&& \hspace*{-5mm} \;=\; \frac{E_{12} \otimes E_{21}}{z-w}-\frac{E_{23}
\otimes E_{32}}{z-w} -\frac{E_{12} \otimes E_{32}}{z-w-\frac{1}{2}}
+\frac{E_{23} \otimes E_{21}}{z-w+\frac{1}{2}} +
\frac{4(z-w)+3}{(z-w)(2(z-w)+1)}E_{13} \otimes E_{31}
\end{eqnarray*}
The explicit form of $\cR_F$ is proven analogously. \\
The calculation of $\cR_{H}$ is standard, but one has to use the following
formula introduced in \cite{KT}:
\begin{equation}
\label{eq:formulahorribilis}
\prod_{n \ge 0} \exp \left\{ \sum_{i \ge 0} \Big( \frac{1}{u-\gamma} \,
\Big)_{i} \Big( \ln \; \frac{x-\alpha+Nn+1}{x-\beta+Nn+1} \Big)_{-i-1}
\right\} =
\frac{\Gamma(\frac{\gamma-\beta+1}{N})}{\Gamma(\frac{\gamma-\alpha+1}{N})}
\end{equation}
where $(\psi(u))_{i} = \psi_{i}$ is defined as in theorem \ref{theor:RH}.
\finproof

\bigskip

\textbf{Acknowledgements:} We warmfully thank J.~Avan and A.~Molev for
discussions and advices. Some preliminary computations were done using
the symbolic manipulation program \textsc{Form}, by
J.~Vermaseren~\cite{Form}. 


\end{document}